\documentstyle[sprocl]{article}

\input{epsf}
\def \beq{\begin{equation}}
\def \eq{\end{equation}}
\def \berr{\begin{eqnarray}}
\def \err{\end{eqnarray}}

\def \a{\alpha}

\def \d{\delta}

\def \w{\omega}

\def \l{\lambda}

\def \U{{\cal U}}

\def \R{{\cal R}}

\def \({\left(}
\def \){\right)}
\def \<{\langle}
\def \>{\rangle}
\def \[{\left[}
\def \]{\right]}
\def \obar{\overline}

\def\tens{\mathop{\otimes}}

\def\ttens{\tilde{\tens}}

\newcommand \compl{C \! \! \! \! {\scriptscriptstyle {}^{{}_|}}\ }

\def\reps{representations }
\def\rep{representation }

\newtheorem{prop}{Proposition}[section]
\newtheorem{theorem}[prop]{Theorem}

\begin{document}

    October 1997
    \hfill  LMU - TPW 97 - 25  
    \vskip .4in
  
\title{Unitary Representations and BRST Structure of the Quantum 
  Anti--de Sitter Group  at Roots of Unity}

\author{ Harold Steinacker }

\address{Sektion Physik der LMU M\"unchen, LS Prof. Wess \\
Theresienstr. 37, D-80333 Munich, Germany}

\maketitle\abstracts{It is shown that for suitable roots of unity, there exist
finite--dimensional unitary \reps of $U_q(so(2,3))$ corresponding to all 
classical one--particle \reps with (half)integer spin, with the correct low--energy
limit. In the massless case for spin $\geq 1$,  a subspace of ``pure gauges'' appears 
which must be factored out, as classically. Unitary many--particle \reps
are defined, with the same low--energy states  as classically.  Furthermore, a 
remarkable element of the center of $U_q(so(2,3))$ is identified which 
plays the role
of the BRST operator, for any spin. The corresponding ghosts are an 
intrinsic part 
of indecomposable representations. }
  
\section{One--Particle  Representations of $U_q(so(2,3))$}

The quantum Anti--de Sitter (AdS) group $\U:=U_q(so(2,3))$ 
is defined as a real form of the 
Drinfeld--Jimbo quantized universal enveloping algebra 
$U_q(so(5,\compl))$, which  
is the Hopf  algebra
\cite{jimbo,drinfeld}
\beq
\[H_i, H_j\] = 0, \quad  \[H_i, X^{\pm}_j\] =
 \pm A_{ji} X^{\pm}_j, \quad \[X^+_i, X^-_j\]   = \d_{i,j} [H_i]_{q_i}    
\label{UEA}
\eq
for $i=1,2$, plus the quantum Serre relations,
where $A_{ij}$ is the Cartan matrix of $B_2$, 
$q_i = q^{\frac{(\a_i,\a_i)}2}$, and 
$[n]_{q_i}  = \frac{q_i^n-q_i^{-n}}{q_i-q_i^{-1}}$.
Coproduct, antipode and counit are defined as usual.
Root vectors corresponding to $\a_3=\a_1+\a_2, \; 
\a_4=\a_1+2\a_2$ can be defined as well.
The Drinfeld--Casimir $v= (S\R_2) \R_1 q^{-2\rho}$ 
is constructed using the universal $\R \in \U \tens \U$,
where $\rho=\frac 32 H_1 + H_2$.

Since we are interested in the root of unity case,  the reality structure of
$\U$ for $|q|=1$ is chosen as
$\obar{H_i} = H_i, \quad  \obar{X^+_1} = -X^-_1, \quad \obar{X^+_2}
= X^-_2$, which is compatible with the coproduct if it also flips the factors in
a tensor product \cite{thesis}.
Then $E  = H_1+\frac 12 H_2 $ is the energy, and
$J_z = \frac 12 H_2$ is a component of angular momentum.

Denote the irreducible \reps (irreps) with lowest weight 
$\mu=E_0 \a_3-s \a_2 \equiv (E_0,s)$ as $V_{(\mu)}$. Elementary particles 
on AdS space  are unitary $V_{(\mu)}$
with $E_0 >0$. The most important ones are the $V_{(\mu)}$  with  
integral $\mu$ and $E_0 \geq s+1 \geq 0$ \cite{fronsdal}.

In the classical case, they are of course infinite--dimensional. 
In the massless case $E_0 = s+1$,
the "naive" lowest weight representations are not completely
reducible, they
develop an invariant subspace of "pure gauge" states with lowest weight 
$(E_0+1, s-1)$ \cite{fronsdal}. To obtain unitarizable
irreps (photon, graviton, ...), one has to factor them out.

In the quantum case for $|q|=1$, the corresponding \reps can only be 
unitarizable if $q$ is a root of unity. Thus from now on $q=e^{2\pi i n/m}$.
Then the $V_{(\mu)}$ can be obtained from compact \reps by a "shift":
Let $V(\l)$ be a highest weight irrep
with {\em compact} highest weight   $\l$, i.e.  
$\l$ is integral with $0 \leq (\l,\a_i) <\frac m{2n} $ for $i=1,...,4$.
Then $\w := V_{(\l_0)}$ with $\l_0=\frac m{2n} \a_3$ is one--dimensional, and
$V_{(\mu)}:=V(\l) \tens \w$  is a positive energy irrep with
lowest weight $\mu =  -\l +\l_0 $.

We call $V_{(\mu)}$ {\em physical} if it is unitarizable w.r.t. $U_q(so(2,3))$.
If $n=1$,  $V_{(\mu)}$ is called {\em Di} for $\mu=(1,1/2)$, and {\em Rac} for
$\mu=(1/2,0)$. Now one can show \cite{unitary}

\begin{theorem}
All $V_{(\mu)}$ where  
$-(\mu - \frac m{2n} \a_3)$ is compact
are physical, in particular the massless irreps, as well as
the singletons Di and Rac. For $E \leq \frac m{2n}$,
they are obtained from a (lowest weight)
Verma module by factoring out the 
submodule with lowest weight $(E_0, -(s+1))$ only, except for the massless 
case, where an additional lowest weight state 
with weight $(E_0+1, s-1)$ appears, and for the Di resp. Rac, where an 
additional lowest weight state with weight $(E_0+1,s)$ resp. $(E_0+2,s)$ 
appears.
This is the same as classically.
\label{noncpct_thm}
\end{theorem} 
For the singletons, this was already shown in \cite{dobrev}.
All these \reps are now finite--dimensional
(we essentially work in the unrestricted specialization).


\section{Tensor Product and Many--Particle Representations}
\label{sec:many_4d}

In general, the full tensor product of 2 physical irreps is not unitarizable,
and indecomposable at 
roots of unity. The idea is to keep the appropriate  
physical lowest weight quotient modules only. 
Consider two physical irreps 
$V_{(\mu)}$ and $V_{(\mu')}$. 
For a basis $\{u_{\l'}\}$ of the physical lowest 
weight states in $V_{(\mu)} \tens V_{(\mu')}$, 
let $Q_{\mu,\mu'}$ be the quotient of $\sum \U \cdot u_{\l'}$ 
after factoring out 
all submodules of the $\U \cdot u_{\l'}$. 
Then $Q_{\mu,\mu'} = \oplus V_{(\l'')}$ where $V_{(\l'')}$ are 
physical lowest weight irreps. Therefore we can define 
\beq
V_{(\mu)} \ttens V_{(\mu')} := \bigoplus_{\l''} V_{(\l'')}.
\eq
If $\frac m{2n}$ is not integer, then the physical states
have non--integral weights, and  
$V_{(\mu)} \ttens V_{(\mu')}$ is zero.  Now we have

\begin{theorem} 
If all weights $\mu, \mu', \dots$ involved are integral, then
$\ttens$ is associative, and 
$V_{(\mu)} \ttens V_{(\mu')}$ is unitarizable w.r.t. $U_q(so(2,3))$.
\end{theorem}

In particular for  $q=e^{2\pi i/m}$ with $m$ even, none of the low--energy
states have been discarded, and our definition is physically sensible.

\section{BRST structure}

To describe spin 1 particles,
consider $V_{(\l)} \tens V_5$, where $V_5$ is the 5 dimensional \rep
of $\U$ ("basic one--forms"), for $q$ as above. In the massive case, this 
is completely reducible as classically, see figure 1a). In the massless case
however,  the would--be irrep with lowest weight $\l-\a_2$ 
has an invariant subspace as in Theorem \ref{noncpct_thm}, which in fact
combines with $V_{(\l+\a_3)}$ 
into one reducible, but indecomposable 
representation, with structure as in figure 1 b).

\begin{figure}
\begin{center}
\leavevmode
\epsfysize=1.9in 
\epsfbox{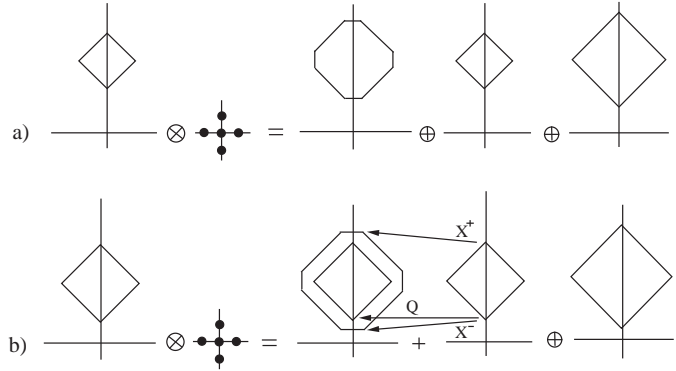}
\caption{$V_{(\l)} \tens V_5$ for the a) massive and 
b) massless case. The $+$
in b) means sum as vector spaces, but not as representations.}
\end{center}
\label{fig:ads_tens_q}
\end{figure}

The main observation \cite{thesis}  now is that this "pure gauge" subspace
is precisely the image of a "BRST" operator 
\beq
Q := (v^{2m} - v^{-2m}),
\label{brst}
\eq
acting on $V_{(\l)} \tens V_5$ . $Q$ vanishes on
any irrep, and $Q^2=0$ on  $V_{(\l)} \tens V_5$ . Thus the physical
Hibert space can be defined as $\mbox{Ker(Q)}/\mbox{Im(Q)}$.
The same $Q$ works similarly for massless particles with higher spin.



\end{document}